\documentclass[12pt]{amsart}
\usepackage{amsmath,amssymb,amsthm,amscd}

\def\p{\partial}
\def\R{\mathbb{R}}
\def\C{\mathbb{C}}
\def\N{\mathbb{N}}

\def\ve{\varepsilon}

\def\L{\Lambda}

\def\i{\sqrt{-1}}

\def\t{\triangle}

\numberwithin{equation}{section}

\newtheorem{prop}{Proposition}[section]
\newtheorem{theo}[prop]{Theorem}

\newtheorem{rmk}[prop]{Remark}

\begin{document}
\title{On  regularity of complex Monge-Amp\`ere equations}
\author{Weiyong HE}

\begin{abstract}
We shall consider the regularity of solutions for complex Monge-Amp\`ere equations in $\C^n$ or a bounded domain. First we prove interior $C^2$ estimates of solutions in a bounded domain for complex Monge-Amp\`ere equations with assumption of  $L^p$ bound for $\t u$, $p>n^2$  and  of Lipschitz condition on right hand side. Then we shall construct  a family of Pogorelov-type examples for complex Monge-Amp\`ere equations. These examples give   generalized entire solutions (as well as  viscosity solutions) of complex Monge-Amp\`ere equation $\det(u_{i\bar j})=1$ in $\C^n$.
\end{abstract}

\address{Department of Mathematics, University of Oregon, Eugene, OR, 97403}
\email{whe@uoregon.edu}

\thanks{2010 MSC: 35J60, 35J96}
\date{}
\maketitle

\section{Introduction}
Let $u$ be a continuous plurisubharmonic function in $\C^n$ (or a bounded domain $\Omega\subset \C^n$). We shall consider  complex Monge-Amp\`ere equations 
\begin{equation}\label{E-1-1}\det\left(\frac{\p^2 u}{\p z_i \p {\bar {z_j}}}\right)=\psi.\end{equation}  For simplicity we use the notations
\[
\begin{split}
u_k=\frac{\p u}{\p z_k},\; u_{\bar k}=\frac{\p u}{\p \bar z_{k}}, \; u_{i\bar j}=\frac{\p^2 u}{\p z_i\p \bar z_{j}}, \; \mbox{etc}.
\end{split}
\]
We use $\t u=\sum_k u_{k\bar k}$ to denote the complex Laplacian operator.

These equations \eqref{E-1-1} have been studied extensively by many mathematicians, for example \cite{Bedford-Taylor1, Bedford-Taylor2, Bedford-Taylor3, Ca-Kohn-Ni-Sp, Guan,
R-Schulz, Schulz} etc. and many others.  In particular, Bedford-Taylor \cite{Bedford-Taylor1}  introduced generalized solutions for complex Monge-Amp\`ere equations;  when $\Omega$ is a strongly pseudo-convex bounded domain, Bedford-Taylor \cite{Bedford-Taylor1, Bedford-Taylor2, Bedford-Taylor3}  established  the existence, uniqueness, and global Lipschitz regularity of generalized solutions for the Dirichlet problem 
\begin{equation}\label{E-1-2}
\begin{split}
& \det(u_{i\bar j})=\psi(z, u, \nabla u), ~\mbox{in} ~\Omega \\
& u=\varphi, ~\mbox{on} ~ \p \Omega.
\end{split}
\end{equation}
Caffarelli-Kohn-Nirenberg-Spruck \cite{Ca-Kohn-Ni-Sp} proved the existence of classical plurisubharmonic solutions of \eqref{E-1-2} under suitable conditions on $\psi$, $\varphi$. 

In the present note we shall first consider the {\it interior}  {\it a priori} estimates of complex Monge-Amp\`ere equations 
\begin{equation}\label{cma-1}
\log \det(u_{i\bar j})=F
\end{equation}
in a bounded domain $\Omega\subset \C^n (n\geq 2)$. The solution $u(z)$,  of class $C^3(\Omega)$,  is strictly plurisubharmonic, such that 
\[
\left(u_{i\bar j}\right) >0, \forall z \in \Omega.
\]
The given function $F$ is of class $C^1(\Omega)$ and satisfies, for some positive constant $\L$,
\[
|F|\leq \L, |\nabla F|\leq \L.
\]
We can state our interior estimates as follows
\begin{theo}\label{T-1} For any $P_0>n^2$ and any domain $\Omega^{'}\subset\subset \Omega$, if $\|\t u\|_{L^{P_0}(\Omega)}$ is bounded, then 
\[
\sup_{\Omega^{'}}\t u \leq C_1=C_1(\L, \Omega, \Omega^{'}, P_0, n, \|\t u\|_{L^{P_0}(\Omega)}).
\]
If we assume that $\|u\|_{W^{2, P_0}(\Omega)}$ is bounded, then it follows that for any $\alpha\in (0, 1)$
\[
\|u\|_{C^{2, \alpha}(\Omega^{'})}\leq C_2=C_2(\L, \Omega, \Omega^{'}, P_0, \alpha, n, \|u\|_{W^{2, P_0}(\Omega)}),
\]
and for any $p>1$,
\[
\|u\|_{W^{3, p}(\Omega^{'})}\leq C_2=C_2(\L, \Omega, \Omega^{'}, P_0, p, n, \|u\|_{W^{2, P_0}(\Omega)}).
\]
\end{theo}

\begin{rmk}
Z. Blocki \cite{Blocki99} proved a similar interior  regularity theorem by assuming $F\in W^{2, 2n}$ and $u\in W^{2, p}$ for $p>2n(n-1)$.  The interior regularity does not hold only assuming $P_0< n(n-1)$ even when $F$ is smooth (see the following example); one might expect that $P_0>n(n-1)$ is the optimal choice for the interior $C^2$ regulairty by considering the following example of  Pogorelov's type, where $u$ in $B_\epsilon(0)$ for $\epsilon$ small enough, 
\[
u(z_1, \cdots, z_n)=(1+|z_1|^2)(|z_2|^2+\cdots+|z_n|^2)^{1-1/n}.
\]
This example was first considered in \cite{Blocki99} by Z. Blocki.
\end{rmk}

For fully nonlinear equations, the regularity of right hand side $F$ does not imply necessarily the regularity of weak solutions. Pogorelov  constructed a well-known example to the real Monge-Amp\`ere equation with $u(x)=(1+x_1^2)(x_2^2+\cdots+x_n^2)^{1-1/n}$ for $n\geq 3$ in $B_{\epsilon}(0)$ for $\epsilon$ small; $u(x)$ is a weak solution of $D^2u=f$ in $B_\epsilon(0)$ such that $0<f\in C^\infty$, while  $u\in C^{1, \alpha}$  for $\alpha=1-2/n$, but not any larger $\alpha$. However, the interior $C^2$ estimates can be established if one assumes additional hypotheses, such as sufficient regularity of $\p \Omega$ and $u|_{\p \Omega}$, or sufficient interior regularity.  In \cite{Guan}, B. Guan studied \eqref{E-1-2} for general domains assuming the existence of a sub-solution. In particular he proved, among others, the interior $C^2$ estimates (depending on $\psi$ up to its second derivatives).  When $F$ is assumed to be only Lipschitz and $u|_{\p \Omega}\equiv 0$, the interior $C^2$ estimates were studied by F. Schulz \cite{Schulz} by using integral  approach of N.M. Ivochikina \cite{Ivo} to the real Monge-Amp\`ere equations. However, the proof in \cite{Schulz} is not complete (this was first pointed out by Z. Blocki \cite{Blocki03}), where the inequality between $(6)$ and $(7)$ in \cite{Schulz} does not follow in the complex case. By using the integral approach, we shall prove the interior $C^2$ estimates by assuming $L^{P_0}$ bound of $\t u$ for some $P_0>n^2$. For real Monge-Amp\`ere equations,  J. Urbas proved $C^{2, \beta}$ regularity of weak solution for any $\beta\in (0, 1)$  (see \cite{Urbas}) assuming $u\in C^{1, \alpha}$ or $u\in W^{2, q}$ for $\alpha>1-2/n, q>n(n-1)/2$. Note that  the bound on $\alpha$ in Urbas's result is sharp for $n\geq 3$; also for $q<n(n-1)/2$, the interior $C^2$ estimates fail for $n\geq 3$. 

Then we shall construct  Pogorelov-type examples for complex Monge-Amp\`ere equations.  Using these examples we construct a generalized solution ( Bedford-Taylor \cite{Bedford-Taylor1}), which is also a viscosity solution,  of the complex Monge-Amp\`ere equation in $\C^n$ such that
\begin{equation}\label{E-3-1}
\det(u_{i\bar j})=1.
\end{equation}
 We have
 \begin{theo}\label{T-3-2}
 Let $z=(z_1, \cdots, z_n)\in \C^n$ and
 \[
 v(z_1, \cdots, z_n)=n^{2/n}(1+|z_1|^2+\cdots+|z_{n-1}|^2) |z_{n}|^{2/n}.
 \]Then $v(z)$ is a generalized solution and a viscosity solution of \eqref{E-3-1}.
 \end{theo}
\begin{rmk} The examples in Theorem \ref{T-3-2} were also considered by Z. Blocki \cite{Blocki99} when $n=2$. When $u_{i\bar j}$ is assumed to be bounded,
Riebesehl-Schulz \cite{R-Schulz} proved that $u_{i\bar j}$ are all constants  if a strictly plurisubharmonic function $u$ solves \eqref{E-3-1}. 
 \end{rmk}
\vspace{2mm}

{\bf Acknowledgement}: I am grateful to Prof. Xiuxiong Chen and Prof. Jingyi Chen for constant support and encouragements. I  benefit from conversations with Prof. Pengfei Guan, Song Sun and Prof. Yu Yuan about complex Monge-Amp\`ere equations; I would like to thank all of them. I would also like to thank Prof. Blocki for pointing out the references \cite{Blocki99, Blocki03} to me and for his interest in this paper.  I am partially supported by a startup grant of University of Oregon. 

\section{Interior $C^2$ estimates}
In this section we prove Theorem \ref{T-1}. Note that $u$ is of class $C^3$ and $F$ is of class $C^1$, in the following $(\t u)_{i\bar j},$ $\t F$ etc. will be understood in the sense of distribution.  
\begin{proof}
We compute, by taking derivative of \eqref{cma-1},
\[
u^{i\bar j} u_{i\bar j k}=F_k.
\]
Taking derivatives again, we get 
\begin{equation}\label{L-2-2}
-u^{i\bar b}u^{a\bar j}u_{a\bar b \bar k}u_{i\bar j k}+u^{i\bar j} (\t u)_{i\bar j}=\t F.
\end{equation}
Let $\phi$ be a test function with support in $\Omega$. We compute the following  for some $p\geq 1$, 
\begin{equation}\label{L-2-3}
\int_\Omega\phi^2 u^{i\bar j} ( (\t u)^p)_i ((\t u)^p)_{\bar j} \t u =p^2 \int_\Omega \phi^2 (\t u)^{2p-1}u^{i\bar j} (\t u)_i (\t u)_{\bar j}.
\end{equation}
Integration by parts, we can get that
\begin{equation}\label{L-2-4}
\begin{split}
&\int_\Omega \phi^2 (\t u)^{2p-1}u^{i\bar j} (\t u)_i (\t u)_{\bar j}\\&=-2\int_\Omega \phi \phi_{\bar j} (\t u)^{2p}u^{i\bar j}  (\t u)_i
-(2p-1)\int_\Omega \phi^2 (\t u)^{2p-1} u^{i\bar j} (\t u)_i (\t u)_{\bar j}\\
&\quad+\int_\Omega \phi^2 (\t u)^{2p} u^{i\bar l} F_{\bar l} (\t u)_i-\int_\Omega \phi^2 (\t u)^{2p} u^{i\bar j} (\t u)_{i\bar j}.
\end{split}
\end{equation}
It follows from \eqref{L-2-2} and \eqref{L-2-4} that 
\begin{equation}\label{L-2-5}
\begin{split}
&2p \int_\Omega \phi^2 (\t u)^{2p-1}u^{i\bar j} (\t u)_i (\t u)_{\bar j}\\&=-2 \int_\Omega \phi \phi_{\bar j} (\t u)^{2p}u^{i\bar j}  (\t u)_i
+\int_\Omega \phi^2 (\t u)^{2p} u^{i\bar l} f_{\bar l} (\t u)_i\\
&\quad -\int_\Omega \phi^2 (\t u)^{2p} \t F-\int_\Omega \phi^2 (\t u)^{2p} u^{i\bar b}u^{a\bar j}u_{a\bar b \bar k}u_{i\bar j k}.
\end{split}
\end{equation}
Recall that $\t F$ is understood in distribution sense. Hence we have
\[
\begin{split}
-\int_\Omega \phi^2 (\t u)^{2p} \t F=& \int_\Omega \left(\phi^2 (\t u)^{2p} \right)_{i} F_{\bar i}\\
=& \int_\Omega 2\phi \phi_{i} (\t u)^{2p} F_{\bar i}+\int_\Omega 2p \phi^2 (\t u)^{2p-1} (\t u)_i F_{\bar i} 
\end{split}
\]
Insert this into \eqref{L-2-5}, it follows that

\begin{equation}\label{L-2-6}
\begin{split}
&2p \int_\Omega \phi^2 (\t u)^{2p-1}u^{i\bar j} (\t u)_i (\t u)_{\bar j}\\&=-2 \int_\Omega \phi \phi_{\bar j} (\t u)^{2p}u^{i\bar j}  (\t u)_i
+\int_\Omega \phi^2 (\t u)^{2p} u^{i\bar l} f_{\bar l} (\t u)_i+\int_\Omega 2\phi \phi_{i} (\t u)^{2p} F_{\bar i}\\
&\quad+\int_\Omega 2p \phi^2 (\t u)^{2p-1} (\t u)_i F_{\bar i} 
 -\int_\Omega \phi^2 (\t u)^{2p} u^{i\bar b}u^{a\bar j}u_{a\bar b \bar k}u_{i\bar j k}.
\end{split}
\end{equation}

To estimate the right hand side in \eqref{L-2-6}, we choose a coordinate such that at one point, $u_{i\bar j}$ is diagonalized. 
Then we compute, at the point, 
\[
\begin{split}
\left|\phi \phi_{\bar j} (\t u)^{2p} u^{i\bar j} (\t u)_i\right|&=\left|\sum_i\frac{1}{u_{i\bar i}}\phi\phi_{\bar i}(\t u)^{2p}(\t u)_i\right|\\
&=\left|\sum_i\frac{1}{u_{i\bar i}}\phi \phi_{\bar i}(\t u)^{2p}\sum_k u_{i\bar k k}\right|\\
&\leq \ve \phi^2 (\t u)^{2p} \sum_{i, k} \frac{|u_{i\bar k k}|^2}{u_{i\bar i}u_{k\bar k}}+\frac{1}{4\ve} |\nabla \phi|^2(\t u)^{2p}\sum_{i, k} \frac{u_{k\bar k}}{u_{i\bar i}}\\
&\leq \ve \phi^2(\t u)^{2p}u^{i\bar b}u^{a\bar j}u_{a\bar b \bar k}u_{i\bar j k}+\frac{e^{-\L}}{4\ve}|\nabla \phi|^2 (\t u)^{2p+n}
\end{split}
\]
Similarly we compute,
\[
\left|\phi^2 (\t u)^{2p} u^{i\bar l} f_{\bar l} (\t u)_i\right|\leq \ve \phi^2(\t u)^{2p}u^{i\bar b}u^{a\bar j}u_{a\bar b \bar k}u_{i\bar j k}+\frac{e^{-\L}}{4\ve}\phi^2 |\nabla F|^2(\t u)^{2p+n},
\]
and
\[
\left|2p \phi^2 (\t u)^{2p-1} (\t u)_i F_{\bar i}\right|\leq \ve \phi^2(\t u)^{2p}u^{i\bar b}u^{a\bar j}u_{a\bar b \bar k}u_{i\bar j k}+e^{-\L}\frac{p^2}{\ve} \phi^2(\t u)^{2p+n-2} |\nabla F|^2.
\]

Insert the above into \eqref{L-2-6} (for example we can choose $\ve=1/4$ in the above inequalities), we can get
\[
2p \int_\Omega \phi^2 (\t u)^{2p-1}u^{i\bar j} (\t u)_i (\t u)_{\bar j}\leq C\int_\Omega(\phi^2+|\nabla \phi|^2) (\t u)^{2p+n}+C p^2\int_\Omega \phi^2(\t u)^{2p+n-2},
\]
where $C$ denotes a universal constant depending only on $\L, n$, and it can vary line by line.  Insert the above into \eqref{L-2-3}, we can get that
\begin{equation}\label{L-2-7}
\int_\Omega\phi^2 u^{i\bar j} ( (\t u)^p)_i ((\t u)^p)_{\bar j} \t u\leq C p^3\int_\Omega(\phi^2+|\nabla \phi|^2) (\t u)^{2p+n}.
\end{equation}
Note that $u^{i\bar j}\t u \xi_i \bar \xi_{j}\geq |\xi|^2$ for any vector $\xi=(\xi_1, \cdots, \xi_n)$. We can get that 
\begin{equation}\label{L-2-8}
\int_\Omega \phi^2\left |\nabla\left((\t u)^p\right)\right|^2\leq C p^3\int_\Omega(\phi^2+|\nabla \phi|^2) (\t u)^{2p+n}.
\end{equation}
It follows that, from \eqref{L-2-8},
\begin{equation}\label{L-2-9}
\int_\Omega \left |\nabla\left(\phi (\t u)^p\right)\right|^2\leq C p^3\int_\Omega(\phi^2+|\nabla \phi|^2) (\t u)^{2p+n},
\end{equation}
where $C$ depends only on $\L, n$.
Starting from \eqref{L-2-9}, one can use the iteration technique (Moser iteration) to get interior $L^\infty$ bound for $\t u$.  By Sobolev inequality, there exists a constant $c=c(n)$ such that 
\begin{equation*}
\|\phi (\t u)^p\|_{L^{2n/(n-1)}(\Omega)}\leq c \|\nabla(\phi (\t u)^p)\|_{L^2(\Omega)}.
\end{equation*}
It follows that
\begin{equation}\label{L-2-10}
\|\phi (\t u)^p\|_{L^{2n/(n-1)}(\Omega)}\leq Cp^{3/2} \left(\int_\Omega(\phi^2+|\nabla \phi|^2) (\t u)^{2p+n}\right)^{1/2}.
\end{equation}
First we consider $\Omega=B_R(0), \Omega^{'}=B_r(0)$ for $0<r<R\leq 1$.
Now we can choose the cut-off function $\phi$ and $p$. Let $r_k=r+(R-r)2^{-k}$, $k\in \N$. Let $p=p_k$ and $\phi=\phi_k$ in \eqref{L-2-10} such that
$\phi_k\equiv 1$ in $B_{r_k}(0)$, $\phi_k\equiv 0$ outside $B_{r_{k-1}}(0)$, and $|\nabla \phi_k|\leq 2(r_{k-1}-r_k)^{-1}$.
By \eqref{L-2-10}, we get
\begin{equation}\label{L-2-11}
\|\t u\|_{L^{\frac{2np_k}{n-1}}(B_{r_k}(0))}\leq a_k \|\t u\|_{L^{2p_k+n}(B_{r_{k-1}(0)})}^{b_k},
\end{equation}
where 
\[a_k=\left(\frac{C p_k^{3/2}}{r_{k-1}-r_k}\right)^{1/p_k}=\left(\frac{C 2^kp_k^{3/2}}{R-r}\right)^{1/p_k}, ~~b_k=\frac{2p_k+n}{2p_k}.\]
Now let $p_0=a+n(n-1)/2, a>0$ and $p_k, k\in \N$ satisfy 
\[
2p_{k+1}+n=\frac{2np_k}{n-1}.
\]
It is easy to get that
\[
p_k=a \left(\frac{n}{n-1}\right)^k+\frac{n(n-1)}{2}.
\]
By \eqref{L-2-11}, we can get that
\begin{equation}\label{L-2-12}
\left\|\t u\right\|_{L^{\frac{2np_k}{n-1}}(B_{r_k}(0))}\leq a_k a_{k-1}^{b_k} a_{k-2}^{b_kb_{k-1}}\cdots a_1^{b_k\cdots b_2} \left\|\t u\right\|_{L^{\frac{2np_0}{n-1}}(B_{r_0}(0))}^{b_k\cdots b_1}.
\end{equation}
It is easy to show that
\[
\lim_{k\rightarrow \infty}\Pi_{i=1}^{k} b_i =p_0/a.
\]
Without loss of generality, we assume $a_k\geq 1, \forall k$.  
It is also clear that
\[
\sum_k \log a_k \leq C=C(R, r, a, n, \L),
\]
which implies that $a_k a_{k-1}^{b_k} a_{k-2}^{b_kb_{k-1}}\cdots a_1^{b_k\cdots b_2}\leq C=C(R, r, a, n, \L)$. 
In \eqref{L-2-12}, let $k\rightarrow \infty$, we can get 
\begin{equation}\label{L-2-13}
\|\t u\|_{L^\infty(B_r(0))}\leq C(R, r, a, n, \L) \|\t u\|^{p_0/a}_{L^{P_0}(B_R(0))},
\end{equation}
where $P_0=2np_0/(n-1)=n^2+2na/(n-1)$. For general domains $\Omega^{'}\subset\subset \Omega$, we can use finite many concentric balls  $b_i, B_i$ such that $b_i\subset\subset B_i\subset\subset \Omega$ and $\overline{\Omega^{'}}\subset \cup b_i$.   Apply \eqref{L-2-13}, one can easily prove that
\[
\|\t u\|_{L^\infty(\Omega^{'})}\leq C_1=C_1(\L, \Omega, \Omega^{'}, P_0, n, \|\t u\|_{L^{P_0}(\Omega)}).
\]
If we assume that $\|u\|_{W^{2, P_0}(\Omega)}$ is bounded, then after obtaining interior $C^2$ estimates as above, one can prove interior $C^{2, \alpha}$ or $W^{3, p}$ estimates for any $\alpha\in (0, 1)$ or $p>1$ by using  Evans-Krylov theory, $L^p$ theory and Schauder theory, for example see \cite{Gil-Tru}.
\end{proof}

\section{Pogorelov-type examples}
In this section we prove Theorem \ref{T-3-2}.
 \begin{proof}
 First we consider $n=2$. Let \begin{equation}u^\epsilon(z, w)=2(1+|z|^2)(|w|^2+\epsilon)^{1/2}, \epsilon>0, (z, w)\in C^2 \end{equation} be a family of smooth functions defined in $\C^2$. A straightforward computation yields
\[\begin{split}
u^{\epsilon}_{z\bar z}&=2(|w|^2+\epsilon)^{1/2}, u^\epsilon_{z\bar w}=\bar z w(|w|^2+\epsilon)^{-1/2},\\
 u_{w\bar w}^\epsilon&=\frac{1+|z|^2}{2}(|w|^2+\epsilon)^{-3/2}(|w|^2+2\epsilon).
\end{split}
\]
It is clear that $u^\epsilon (\epsilon>0)$ is a smooth plurisubharmonic function defined in $\C^2$ and satisfies the complex Monge-Amp\`ere equation
\begin{equation}
u^\epsilon_{z\bar z} u^\epsilon_{w\bar w}-u^{\epsilon}_{z\bar w}u^{\epsilon}_{w\bar z}=F(z, w, \epsilon),
\end{equation}
where
\[F(z, w, \epsilon)=|w|^2(|w|^2+\epsilon)^{-1}+2\epsilon(1+|z|^2)(|w|^2+\epsilon)^{-1}.\] 
When $\epsilon\rightarrow 0$, $u^\epsilon$ converges to a continuous plurisubharmonic function $u^0=(1+|z|^2)|w|$. It is clear that $u^0\in Lip_{loc}(\C^2)\cap W^{2, p}_{loc}(\C^2)$ for any $p\in (0, 2)$. So  $u^\epsilon$ converges to $u^0$ in $C^\alpha$ and $W^{2, p}$ norm on any compact subset of $\C^2$  for any $\alpha\in (0, 1)$ and $p\in (0, 2)$. 
It is  clear that  $F(z, w, \epsilon)$ converges to $1$ in $L^\infty$ sense when $\epsilon\rightarrow 0$. 
In particular, $u^0$ is a generalized solution  (in the sense of \cite{Bedford-Taylor1})  of complex Monge-Amp\`ere equation in $\C^2$ such that
\begin{equation}\label{E-3-5}
u^0_{z\bar z} u^0_{w\bar w}-u^0_{z\bar w}u^0_{w\bar z}=1.
\end{equation}
Now we show that $u^0$ is also a viscosity solution. We refer to \cite{Cra-Ishii-Lions} for more details about viscosity solutions.
For a $n\times n$ matrix $X$, define
\begin{equation*}
G(X)=\left\{%
\begin{array}{ll}
 &1-\det(X), X\geq 0 \\
 & +\infty, ~~\mbox{otherwise}.
\end{array}%
\right.
\end{equation*}
Let $\p\bar \p u$ denote the complex Hessian of $u$. 
To show that $u^0$ is a viscosity solution of $G(\p\bar \p u)=0$, we need to show that for any $v\in C^2$, if $u-v$ has a local minimum at $p$, then $G(\p\bar \p v)\geq 0$ at $p$; and if $u-v$ has a local maximum at $p$, then $G(\p\bar\p v)\leq 0$ at $p$. Since $u^0$ solves $G(\p\bar \p u)=0$ when $|w|\ne 0$ in classical sense, we just need to check for $|w|=0$. 

Let $z=x_1+\i y_1$, $w=x_2+\i y_2$.  Suppose $u-v$ has a local minimum at $(a, b, 0, 0)$, then 
\[
u^0(x_1, y_1, x_2, y_2)-v(x_1, y_1, x_2, y_2)\geq u(a, b, 0, 0)-v(a, b, 0, 0)
\]
for any $(x_1, y_1, x_2, y_2)$ in a small neighborhood of $(a, b, 0, 0)$. So we get 
\[
v(x_1, y_1, x_2, y_2)-v(a, b, 0, 0)\leq (1+x_1^2+y_1^2)\sqrt{x_2^2+y_2^2}. 
\]
Take $x_2=y_2=0$, it gives that
\[
v(x_1, y_1, 0, 0)-v(a, b, 0, 0)\leq 0
\]
for any $(x_1, y_1)$ in a small neighborhood of $(a, b)$. This implies that $D^2_{x_1, y_1}v(a, b, 0, 0)\leq 0$. In particular 
$v_{z\bar z}(a, b, 0, 0)\leq 0$. Let $X$ be the matrix of complex Hessian of $v$ at $(a, b, 0, 0)$. Either $X\geq 0$, then $X\equiv 0$
and  $G(X)=1>0$; or  $X\geq 0$ does not hold , then $G(X)=\infty>0.$ 

Suppose $u-v$ has a local maximum at $(a, b, 0, 0)$, then
\[
u^0(x_1, y_1, x_2, y_2)-v(x_1, y_1, x_2, y_2)\leq u(a, b, 0, 0)-v(a, b, 0, 0)
\]
So we get 
\[
v(x_1, y_1, x_2, y_2)-v(a, b, 0, 0)\geq (1+x_1^2+y_1^2)\sqrt{x_2^2+y_2^2}. 
\]
for any $(x_1, y_1, x_2, y_2)$ in a small neighborhood of $(a, b, 0, 0)$. Take $x_1=a, y_1=b$, we get that
\[
v(a, b, x_2, y_2)-v(a, b, 0, 0)\geq (1+a^2+b^2)\sqrt{x_2^2+y_2^2}.
\]
But this contradicts $v\in C^2$. So  there is no $C^2$ function $v$ such that $u-v$ has a local maximum.\\

Then we shall consider $\C^{n+1}$ $(n\geq 2)$. For $(z_1, z_2, \cdots, z_n, w)\in \C^{n+1}$, let 
\begin{equation}
u^{\epsilon}(z_1, \cdots, z_n, w)=(1+ |z_1|^2+\cdots +|z_{n}|^2) (|w|^2+\epsilon)^{1/{n}}.
\end{equation}
A straightforward computation gives that
\[
\begin{split}
u^\epsilon_{z_i\bar z_{j}}&=(|w|^2+\epsilon)^{1/(n+1)}\delta_{ij}, u^\epsilon_{z_i\bar w}=\frac{\bar z_i w}{n+1} (|w|^2+\epsilon)^{-n/(n+1)},\\
u^\epsilon_{w\bar w}&=\frac{1}{(n+1)^2}(1+|z_1|^2+\cdots+|z_n|^2)(|w|^2+\epsilon)^{-(2n+1)/(n+1)}(|w|^2+(n+1)\epsilon).
\end{split}
\]
It is clear that $u^\epsilon (\epsilon>0)$ is a smooth plurisubharmonic function  in $\C^{n+1}$ which solves ($z_{n+1}=w$)
\begin{equation}
\det( u^\epsilon_{i\bar j})=F(z_1, \cdots, z_n, w, \epsilon),
\end{equation}
where
\[
F(z_1, \cdots, z_n, w, \epsilon)=\frac{|w|^2}{(n+1)^2}(|w|^2+\epsilon)^{-1}+\frac{\epsilon}{n+1} (1+|z_1|^2+\cdots+|z_n|^2)(|w|^2+\epsilon)^{-1}.
\]
When $\epsilon\rightarrow 0$, $u^\epsilon\rightarrow u^0=(1+|z_1|^2+\cdots+|z_n|^2)|w|^{2/(n+1)}$.
Note that $u^0\in C^\alpha_{loc}\cap W^{2, p}_{loc}$ for $\alpha=2/(n+1)$ and $p<1+1/n$. So $u^\epsilon$ converges to $u^0$ in $C^\alpha$ and $W^{2, p}$ norm for $\alpha=2/(n+1)$ and $p<1+1/n$ on any compact subset of $\C^{n+1}$; while $F(z_1, \cdots, z_n, w, \epsilon)$ converges to $1$ in $L^\infty$.
In particular, $u^0$ is a generalized solution of the complex Monge-Amp\`ere equation in $\C^{n+1}$,
\begin{equation}\label{E-3-8}
\det (u^0_{i\bar j})=1/(n+1)^2.
\end{equation}
It is similar to check that $u^0(z_1, \cdots, z_n, w)$ is a viscosity solution as in $n=2$. Theorem \ref{T-3-2} then follows from \eqref{E-3-5} and \eqref{E-3-8}.
\end{proof}

\begin{rmk} It is clear that $u=n^{2/n}(\sum_{i=1}^{n-1} |z_i|^2)|z_n|^{2/n}$   is a generalized solution of  the degenerated complex Monge-Amp\`ere equation
\[
\det(u_{i\bar j})=0.
\]
\end{rmk}


\begin{thebibliography}{s}
\bibitem{Bedford-Taylor1} E. Bedford, B.A. Taylor;  {\it The Dirichlet problem for a complex Monge-Amp\`re equation.} Invent. Math. 37 (1976), no. 1, 1--44.
\bibitem{Bedford-Taylor2} E. Bedford, B.A. Taylor;  {\it Variational properties of the complex Monge-Amp\`ere equation. I. Dirichlet principle.} Duke Math. J. 45 (1978), no. 2, 375--403.
\bibitem{Bedford-Taylor3} E. Bedford, B.A. Taylor; {\it  Variational properties of the complex Monge-Amp\`ere equation. II. Intrinsic norms.}  Amer. J. Math. 101 (1979), no. 5, 1131--1166.
\bibitem{Blocki99}Z. Blocki; {\it On the regularity of the complex Monge-Amp\`ere operator}, Contemporary Mathematics 222, Complex Geometric Analysis in Pohang, ed. K.-T.Kim, S.G.Krantz, pp.181-189, Amer. Math. Soc. 1999. 
\bibitem{Blocki03} Z. Blocki; {\it Interior regularity of the degenerate Monge-Amp\`ere equation,} Bulletin of the Australian Mathematical Society 68 (2003), 81-92. 
\bibitem{Caf}L. Caffarelli;  {\it Interior $W^{2,p}$ estimates for solutions of the Monge-Ampre equation.} Ann. of Math. (2) 131 (1990), no. 1, 135--150.  
\bibitem{Ca-Kohn-Ni-Sp}L. Caffarelli, J.J. Kohn, L. Nirenberg, J. Spruck; {\it The Dirichlet problem for nonlinear second-order elliptic equations. II. Complex Monge-Amp\`ere, and uniformly elliptic, equations.} Comm. Pure Appl. Math. 38 (1985), no. 2, 209--252.
\bibitem{Cra-Ishii-Lions} M.G. Crandall, H. Hitoshi, P.L. Lions;  {\it User's guide to viscosity solutions of second order partial differential equations.} Bull. Amer. Math. Soc. (N.S.) 27 (1992), no. 1, 1--67.
\bibitem{Gil-Tru} D. Gilbarg, N. Trudinger, {\it Elliptic partial differential equations of second order}, Springer, 1998 Edition.
\bibitem{Guan}B. Guan; {\it The Dirichlet problem for complex Monge-Amp\`ere equations and regularity of the pluri-complex Green function.} Comm. Anal. Geom. 6 (1998), no. 4, 687--703.
\bibitem{Ivo}N.M. Ivochikina; {\it Construction of a priori bounds for convex solutions of the Monge-Amp\`ere equation by integral methods,} Ukrian. Math. J. 30 (1978), 32-38.
\bibitem{R-Schulz} D. Riebesehl, F. Schulz; {\it A priori estimates and a Liouville theorem for complex Mong-Amp\`ere equations.} Math. Z. 186 (1984), no. 1, 57--66. 
\bibitem{Schulz} F. Schulz; {\it A $C^2$ estimate for solutions of complex Monge-Amp\`ere equations},  J. Reine Angew. Math.  348  (1984), 88--93.
\bibitem{Urbas} J. Urbas; {\it Regularity of generalized solutions of Monge-Amp\`ere equations.} Math. Z. 197 (1988), no. 3, 365--393.
\end{thebibliography}
\end{document}